\documentstyle[12pt,intlim]{amsart}

\def\txt#1{\quad\text{#1}\quad}

\def\sprod#1#2{\left\langle#1,#2\right\rangle}
\def\iprod#1#2{\left(#1,#2\right)}
\def\a{\boldsymbol\alpha}
\def\b{\boldsymbol\beta}
\def\d{\boldsymbol\delta}
\def\g{\boldsymbol\gamma}
\def\m{\boldsymbol\mu}
\def\r{\boldsymbol\varrho}
\def\C{\Bbb C}
\def\K{\Bbb K}
\def\N{\Bbb N}
\def\R{\Bbb R}
\def\Z{\Bbb Z}
\newmathalphabet*{\eusm}{eus}{m}{n}

\theoremstyle{plain}
\newtheorem{corollary}{Corollary}	
\newtheorem{proposition}{Proposition}[section]
\newtheorem{theorem}{Theorem}		
\newtheorem{lemma}{Lemma}[section]

\theoremstyle{definition}
\newtheorem{definition}{Definition}	

\theoremstyle{remark}
\newtheorem{remark}{Remark}		

\title{Ideal norms and\\trigonometric orthonormal systems}
\author[J.~Wenzel]{J\"org Wenzel}
\thanks{This article originated from the author's Ph.D.~thesis at the 
University of Jena written under the supervision of A.~Pietsch}
\address{Mathematical Institute\\FSU Jena\\07740 Jena\\Germany}
\email{\tt wenzel@@minet.uni-jena.de}

\begin{document}
\maketitle

\begin{abstract}
 In this article, we characterize the $UMD$--property of a Banach space $X$ 
 by ideal norms associated with trigonometric orthonormal systems.

 The asymptotic behavior of that numerical parameters can be used to decide 
 whether or not $X$ is a $UMD$--space. Moreover, in the negative case, we 
 obtain a measure that shows how far $X$ is from being a $UMD$--space.

 The main result is, that all described parameters are equivalent also in the 
 quantitative setting.
\end{abstract}


\section{Introduction}
The study of sequences of ideal norms can be used to quantify certain 
properties of linear operators. In most cases the boundedness of a sequence 
of ideal norms for a given operator $T$ describes a well--known property, 
whereas, in the non--bounded case, the growth rate of the sequence describes 
how much the operator $T$ deviates from this property.

One particularly interesting case is if two sequences of ideal norms are 
uniformly equivalent. Then the properties given by these sequences are 
equivalent also in the quantitative setting.

We introduce several sequences of ideal norms related to the trigonometric
orthonormal systems. The boundedness of these sequences for the identity map 
of a Banach space $X$ is equivalent to $X$ being $UMD$.

All of these sequences turn out to be uniformly equivalent. As a corollary 
we get that a Banach space $X$ is a $UMD$--space if and only if there exists 
a constant $c\ge0$ such that, for all $x_1,\dots,x_n\in X$, we have
\[ \left(\frac1\pi\int_{-\pi}^{+\pi}\left\|\sum_{k=1}^nx_k\sin kt\right\|^2
 dt\right)^{1/2}\le c\left(\frac1\pi\int_{-\pi}^{+\pi}\left\|
 \sum_{k=1}^nx_k\cos kt\right\|^2dt\right)^{1/2}
\] or, what turns out to be equivalent,
\[ \left(\frac1\pi\int_{-\pi}^{+\pi}\left\|\sum_{k=1}^nx_k\cos kt\right\|^2
 dt\right)^{1/2}\le c\left(\frac1\pi\int_{-\pi}^{+\pi}\left\|
 \sum_{k=1}^nx_k\sin kt\right\|^2dt\right)^{1/2}.
\]


\section{Ideal norms}

Let $X$ and $Y$ be Banach spaces. Since we deal with the exponential system 
$(\exp(it),\dots,\exp(int))$, most of the results only make sence in the 
complex setting. However, they remain true if the exponential system is 
replaced by its real analogue 
\[ (1,\sqrt2\cos t,\dots,\sqrt2\cos nt,\sqrt2\sin t,\dots,\sqrt2\sin nt).
\]

Let $\boldsymbol{\frak L}$ denote the ideal of all bounded linear operators.

For the theory of ideal norms and operator ideals we refer to the monographs
of Pietsch, \cite{PIE:1} and \cite{PIE:2}. For a more general treatment of 
ideal norms associated with orthonormal systems, we refer to \cite{PIW}.

\begin{definition}
 An {\em ideal norm} $\a$ is a function, which assigns to every operator $T$
 between arbitrary Banach spaces a non--negative number $\a(T)$ such that the
 following conditions are satisfied:
 \begin{align*}
  &\a(S+T)\le\a(S)+\a(T)\txt{for all $S,T\in\boldsymbol{\frak L}(X,Y)$,}\\
  &\a(BTA)\le\|B\|\,\a(T)\,\|A\|\!\!\txt{for all $A\in\boldsymbol{\frak 
  L}(X_0,X)$,
   $T\in\boldsymbol{\frak L}(X,Y)$, $B\in\boldsymbol{\frak L}(Y,Y_0)$,}\\
  &\a(T)=0\txt{implies}T=O.
 \end{align*}
\end{definition}

We write $\a(X)$ instead of $\a(I_X)$, where $I_X$ denotes the identity map 
of the Banach space $X$.

If we additionally assume that $\a(\K)\ge1$, where $\K$ is the scalar field
of the real numbers $\R$ or the complex numbers $\C$, then we have 
$\a(T)\ge\|T\|$ for all operators $T\in\boldsymbol{\frak L}$. The assumption 
above is in particular satisfied by all ideal norms considered in this 
article.

If $\a$ is an ideal norm then its {\em dual ideal norm} $\a'$ is defined by
\[ \a'(T):=\a(T').
\]

The ideal norm $\a$ is said to be {\em injective} if
\[ \a(JT)=\a(T)
\] for all $T\in\boldsymbol{\frak L}(X,Y)$ and any metric injection 
$J\in\boldsymbol{\frak L}(Y,Y_0)$. A metric injection $J$ is a linear map 
such that $\|Jy\|=\|y\|$ for all $y\in Y$.

Let $\a$ be an ideal norm and let $c>0$ be a constant. We write
\[ \a\le c
\] if
\[ \a(X)\le c\txt{for all Banach spaces $X$.}
\] It then follows that for all $T\in \boldsymbol{\frak L}$
\[ \a(T)\le\|T\|\,\a(X)\le c\,\|T\|.
\]

Given ideal norms $\a$, $\b$ and $\g$, we write
\[ \a\le\b\circ\g
\] if
\[ \a(ST)\le\b(S)\,\g(T)\txt{for all $T\in\boldsymbol{\frak L}(X,Y)$ and 
$S\in\boldsymbol{\frak L}(Y,Z)$.}
\]

The following concept is essential for the further considerations.

\begin{definition}
 Two sequences of ideal norms $(\a_n)$ and $(\b_n)$ are said to be
 {\em uniformly equivalent} if there exists a constant $c>0$ such that 
 \[ \frac1c\,\a_n(T)\le\b_n(T)\le c\,\a_n(T)
 \] for all $T\in\boldsymbol{\frak L}$.
\end{definition}


\section{Orthonormal systems}

Given any Banach space $X$ and a measure space $(M,\mu)$, let $L_2^X(M,\mu)$
denote the Banach space of all $\mu$--measurable functions $f:M\to X$ for 
which
\[ \|f|L_2\| := \left( \,\int_M \|f(t)\| ^2 d\mu(t) \right) ^{1/2}
\] is finite.

In the following, let 
\[ \eusm A_n=(a_1,\dots,a_n)\txt{and}\eusm B_n=(b_1,\dots,b_n)
\] be orthonormal systems in some Hilbert space $L_2(M,\mu)$ and $L_2(N,\nu)$,
respectively.

For every orthonormal system $\eusm A_n$, we also consider the complex 
conjugate  orthonormal system $\overline{\eusm A}_n$, which consists of the 
functions $\overline a_1,\dots,\overline a_n\in L_2(M,\mu)$. 

For $x_1,\dots, x_n\in X$, we write
\begin{equation}
 \|(x_k)|\eusm A_n\|:=\left(\,\int_M\left\|\sum_{k=1}^nx_ka_k(s)\right\|^2ds
 \right)^{1/2}.
\end{equation}
This expression yields a norm on the $n$--th Cartesian power of $X$.

\begin{proposition}\label{single elements}
 $\|x_h\|\le\|(x_k)|\eusm A_n\|$ for all $h=1,\dots,n$.
\end{proposition}

\begin{pf}
 By the Parseval equation, we have for all $x'\in X'$
 \[ \sum_{k=1}^n|\sprod{x_k}{x'}|^2=\int_M\left|\sum_{k=1}^n\sprod{x_k}{x'}
  a_k(s)\right|^2d\mu(s)
  \le\int_M\left\|\sum_{k=1}^nx_ka_k(s)\right\|^2d\mu(s)\,\|x'\|^2.
 \] Hence
 \[ \|x_h\|=\sup_{\|x'\|\le1}|\sprod{x_h}{x'}|\le\|(x_k)|\eusm A_n\|.
 \]
\end{pf}

\begin{definition}
 For $T\in \boldsymbol{\frak L}(X,Y)$ and $n\in\N$ the ideal norm $\r(T|\eusm 
 B_n,\eusm A_n)$ is defined as the smallest constant $c\ge0$ such that
 \begin{equation}\label{def r}
  \|(Tx_k)|\eusm B_n\|\le c\,\|(x_k)|\eusm A_n\|
 \end{equation}
 whenever $x_1,\dots,x_n\in X$.

 The ideal norm $\d(T|\eusm B_n,\eusm A_n)$ is defined as the smallest 
 constant $c\ge0$ such that
 \begin{equation}\label{def d}
  \|(T\sprod f{\overline a_k})|\eusm B_n\|\le c\,\|f|L_2\|
 \end{equation}
 whenever $f\in L_2^X(M,\mu)$. Here
 \[ \sprod f{\overline a_k} := \int_M f(s) \overline{a_k(s)} \,d\mu(s)
 \] denotes the $k$--th {\em Fourier coefficient} of $f$ with respect to
 $\eusm A_n$.
\end{definition}

\begin{proposition}
 For any three orthonormal systems $\eusm A_n$, $\eusm B_n$ and $\eusm F_n$, 
 we have
 \begin{eqnarray}\label{r le d}
  \r(\eusm B_n,\eusm A_n)&\le&\d(\eusm B_n,\eusm A_n),\\
  \d(\eusm B_n,\eusm A_n)&\le&\r(\eusm B_n,\eusm F_n)
  \circ\d(\eusm F_n,\eusm A_n),\label{d le rd}\\
  \r(\eusm B_n,\eusm A_n)&\le&\r(\eusm B_n,\eusm F_n)
  \circ\r(\eusm F_n,\eusm A_n).
 \end{eqnarray}
\end{proposition}

\begin{pf}
 The first inequality follows by substituting $f=\sum_{k=1}^nx_ka_k$ in the 
 defining inequality \eqref{def d} of $\d(\eusm B_n,\eusm A_n)$. The other 
 inequalities are trivial.
\end{pf}

The next fact is obvious, as well.
\begin{proposition}
 The ideal norms $\d(\eusm B_n,\eusm A_n)$ are injective.
\end{proposition}

The ideal norms $\d(\eusm B_n,\eusm A_n)$ enjoy the following duality
property.

\begin{proposition}
 $\d(\eusm B_n,\eusm A_n)=\d'(\overline{\eusm A}_n,\overline{\eusm B}_n)$.
\end{proposition}

\begin{pf}
 For $T\in\boldsymbol{\frak L}(X,Y)$ and $g\in L_2^{Y'}(N,\nu)$, we let
 \[ c:=\left(\,\int_M\left\|\sum_{k=1}^nT'\sprod g{b_k}\overline{a_k(s)}
  \right\|^2 d\mu(s)\right)^{1/2}
  =\left\|\left.\left(T'\sprod g{b_k}\right)\right|\overline{\eusm A}_n 
  \right\|.
 \] Given $\varepsilon>0$, by \cite[p. 232]{DIN} there exists $f\in 
 L_2^X(M,\mu)$ such that
 \[ c=\int_M\sprod{f(s)}{\sum_{k=1}^nT'\sprod g{b_k}\overline{a_k(s)}}
  \,d\mu(s)
 \] and
 \[ \|f|L_2\|\le1+\varepsilon.
 \] We now obtain
 \begin{align*}
  c&=\int_M\!\int_N\sum_{k=1}^n\sprod{Tf(s)}{g(t)}\overline{a_k(s)}
  b_k(t)\,d\mu(s)d\nu(t)\\
  &=\int_N\sprod{\sum_{k=1}^nT\sprod f{\overline a_k}b_k(t)}{g(t)}
  d\nu(t)\\
  &\le\left(\,\int_N\left\|\sum_{k=1}^nT\sprod f{\overline a_k}b_k(t)
  \right\|^2d\nu(t)\right)^{1/2}\left(\,\int_N\|g(t)\|^2d\nu(t)
  \right)^{1/2}\\
  &\le(1+\varepsilon)\,\d(T|\eusm B_n,\eusm A_n)
  \,\|g|L_2\|.
 \end{align*}
 Letting $\varepsilon$ tend to $0$ yields
 \[ \left\|\left.\left(T'\sprod g{b_k}\right)\right|\overline{\eusm 
  A}_n\right\|\le
  \d(T|\eusm B_n,\eusm A_n)\,\|g|L_2\|.
 \] This proves that 
 \[ \d(T'|\overline{\eusm A}_n,\overline{\eusm B}_n)
  \le\d(T|\eusm B_n,\eusm A_n).
 \]

 Note that $T''K_X=K_YT$, where $K_X$ and $K_Y$ denote the canonical
 embedding from $X$ into $X''$ and from $Y$ into $Y''$, respectively. Using the
 injectivity of $\d(\eusm B_n,\eusm A_n)$ and $\|K_X\|\le1$,
we finally conclude that
 \begin{align*}
  \d(T|\eusm B_n,\eusm A_n)&=\d(K_YT|\eusm B_n,\eusm A_n)
  =\d(T''K_X|\eusm B_n,\eusm A_n)\\
  &\le\d(T''|\eusm B_n,\eusm A_n)\le\d(T'|\overline{\eusm 
  A}_n,\overline{\eusm B}_n)\le\d(T|\eusm B_n,\eusm A_n).
 \end{align*}
\end{pf}

>From the duality property of the ideal norms $\d(\eusm B_n,\eusm A_n)$ 
and \eqref{d le rd}, we get the following result.

\begin{proposition}\label{cor}
 Let $\eusm A_n$ and $\eusm B_n$ as well as $\eusm F_n$ and
 $\eusm G_n$ be orthonormal systems. Then
 \begin{equation*}
  \d(\eusm B_n,\eusm A_n)\le\r(\eusm B_n,\eusm G_n)\circ
  \d(\eusm G_n,\eusm F_n)\circ\r'(\overline{\eusm A}_n,
  \overline{\eusm F}_n).
 \end{equation*}
\end{proposition}

We denote by $\eusm A_n\otimes\eusm B_n$ the orthonormal system in
$L_2(M\times N,\mu\times \nu)$ consisting of the functions $a_k\otimes b_k:
(s,t)\to a_k(s)b_k(t)$ with $k=1,\dots,n$. Note that
\begin{equation}\label{tensornorm}
 \|(x_k)|\eusm A_n\otimes\eusm B_n\|
 =\|(x_k)|\eusm B_n\otimes\eusm A_n\|
 =\left(\,\int_M\!\!\int_N\left\|\sum_{k=1}^nx_ka_k(s)b_k(t)
 \right\|^2\!\!d\mu(s)\,d\nu(t)\!\right)^{1/2}\!\!.
\end{equation}

The following fact turns out to be very useful to formulate various proofs.
\begin{proposition}\label{prop tensor}
 Let $\eusm F_n=(f_1,\dots,f_n)$ be another orthonormal system in the space
 $L_2(R,\varrho)$. Then
 \begin{equation*}
  \r(\eusm B_n\otimes\eusm F_n,\eusm A_n\otimes\eusm F_n)\le
  \r(\eusm B_n,\eusm A_n).
 \end{equation*}
\end{proposition}

\begin{pf}
 Substituting $(x_kf_k(r))$ with $r\in R$ in the defining inequality
 \eqref{def r} of $\r(T|\eusm B_n,\eusm A_n)$, we obtain
 \[ \int_N\left\|\sum_{k=1}^nTx_kb_k(t)f_k(r)\right\|^2d\nu(t)\le
  \r(T|\eusm B_n,\eusm A_n)^2\int_M\left\|\sum_{k=1}^nx_ka_k(s)f_k(r)
  \right\|^2d\mu(s).
 \] Integration over $r\in R$ and taking square roots yields
 \[ \|(Tx_k)|\eusm B_n\otimes\eusm F_n\|\le\r(T|\eusm B_n,\eusm A_n)\,
  \|(x_k)|\eusm A_n\otimes\eusm F_n\|,
 \] which proves the desired result.
\end{pf}


\section{Trigonometric orthonormal systems}

We write
\[ \begin{array}{rcll}
 e_k(t)&:=&\exp(ikt)&\txt{for $k\in\Z$,}\\
 c_k(t)&:=&\sqrt2\cos kt&\txt{for $k\in\N$,}\\
 s_k(t)&:=&\sqrt2\sin kt&\txt{for $k\in\N$.}
 \end{array}
\]

Note that
\[ \eusm E_n:=(e_1,\dots,e_n),\quad\eusm C_n:=(c_1,\dots,c_n)\txt{and}
 \eusm S_n:=(s_1,\dots,s_n)
\] are orthonormal systems in $L_2(-\pi,+\pi)$ equipped with the scalar 
product
\[ \iprod fg:=\frac1{2\pi}\int_{-\pi}^{+\pi}f(t)\overline{g(t)}\,dt.
\]

Moreover, we have
\[ \|(x_k)| \eusm C_n\| = \hspace{-1pt} \left( \frac2\pi \int_0^\pi \left\|
 \sum_{k=1}^n x_k \cos kt \right\| ^2 \hspace{-2pt} dt \right) ^{1/2}
 \hspace{-15pt} \txt{and} \hspace{-5pt} \|(x_k)| \eusm S_n\| = \hspace{-1pt}
 \left( \frac2\pi \int_0^\pi \left\| \sum_{k=1}^n x_k \sin kt \right\| ^2
 \hspace{-2pt} dt \right) ^{1/2} \hspace{-5pt}.
\]

Note that
\[ \|(x_k)| \eusm E_n\| = \left( \frac1{2\pi} \int_{-\pi}^{+\pi} \left\|
 \sum_{k=1}^n x_k \exp(ikt) \right\| ^2 dt \right) ^{1/2}.
\] Hence the substitution $t\to-t$ yields
\begin{equation}\label{conju}
 \|(x_k)|\overline{\eusm E}_n\|=\|(x_k)|\eusm E_n\|.
\end{equation}


\section{Main result}

We are now ready to state the main result.

\begin{theorem}
 The sequences of the following ideal norms are uniformly equivalent:
 \[\d(\eusm E_n,\eusm E_n),\ \d(\eusm S_n,\eusm C_n),\ 
\d(\eusm C_n,\eusm S_n),\ \r(\eusm S_n,\eusm C_n),\ \r(\eusm C_n,\eusm S_n).\]
\end{theorem}

\begin{definition}
 A Banach space $X$ has the {\em $UMD$--property} if there exists a constant 
 $c\ge0$ such that
 \[ \left\|\left.\sum_{k=0}^n\varepsilon_k\,dM_k\right|L_2\right\|\le c\,
  \left\|\left.\sum_{k=0}^ndM_k\right|L_2\right\|
 \] for all martingales $(M_0,M_1,\dots)$ with values in $X$, all $n\in\N$ and 
 all sequences of signs $(\varepsilon_1,\dots,\varepsilon_n)\in\{\pm1\}^n$; 
 see \cite{Bou83}.
\end{definition}

It is known (see \cite{Bou83}, \cite{BUR:1986}, \cite{DEF:1986}) that a 
Banach space $X$ is a $UMD$--space if the sequence of ideal norms 
$\d(\eusm E_n,\eusm E_n)$ is bounded. Hence we get the following corollary.

\begin{corollary}
 Let $X$ be a Banach space. The following conditions are equivalent:
 \begin{enumerate}
  \item
   $X$ is a $UMD$--space.
  \item
   There exists a constant $c\ge0$ such that, for all $f \in L_2^X (-\pi,
   +\pi)$, we have
   \[ \left( \frac1{2\pi} \int_{-\pi}^{+\pi} \left\| \sum_{k=1}^n
    \sprod{f}{\overline e}_k \exp(ikt) \right\|^2 dt \right) ^{1/2} \le c 
    \left( \frac1{2\pi} \int_{-\pi}^{+\pi} \left\|
f(t) \right\| ^2 dt\right) ^{1/2}.
   \]
  \item
   There exists a constant $c\ge0$ such that, for all $x_1,\dots,x_n\in X$, 
   we have
   \[ \left\| \left. \sum_{k=1}^n x_k \sin kt \right| L_2 \right\| \le c 
    \left\| \left. \sum_{k=1}^n x_k \cos kt \right| L_2 \right\|.
   \]
  \item
   There exists a constant $c\ge0$ such that, for all $x_1,\dots,x_n\in X$, 
   we have
   \[ \left\| \left. \sum_{k=1}^n x_k \cos kt \right| L_2 \right\| \le c 
    \left\| \left. \sum_{k=1}^n x_k \sin kt \right| L_2 \right\|.
   \]
 \end{enumerate}
\end{corollary}


\section{Proof of the main result}

\begin{lemma}\label{lem1}
 $\r(\eusm S_n,\eusm S_n\otimes \eusm C_n)\le\sqrt2.$
\end{lemma}

\begin{pf}
 It follows from
 \[ \sin k(t-s)=\sin kt\cos ks-\cos kt\sin ks
 \] and the translation invariance of the Lebesgue measure that for all
 $s\in\R$ and $x_1,\dots,x_n\in X$
 \begin{align*}
  \|(x_k)|\eusm S_n\|^2&=\frac1\pi\int_{-\pi}^{+\pi}\left\|
  \sum_{k=1}^nx_k\sin k(t-s)\right\|^2dt\\
  &\le\frac2\pi\int_{-\pi}^{+\pi}\left\|\sum_{k=1}^nx_k\sin kt\cos 
  ks\right\|^2dt+
  \frac2\pi\int_{-\pi}^{+\pi}\left\|\sum_{k=1}^nx_k\cos kt\sin ks\right\|^2dt.
 \end{align*}
 Integrating this inequality over $s\in[-\pi,+\pi]$, we get
 \[ \|(x_k)|\eusm S_n\|^2\le\|(x_k)|\eusm S_n\otimes \eusm C_n\|^2+\|(x_k)|
  \eusm C_n\otimes\eusm S_n\|^2.
 \]
 This proves the assertion by taking into account $\|(x_k)|\eusm S_n\otimes
 \eusm C_n\|=\|(x_k)|\eusm C_n\otimes\eusm S_n\|$.
\end{pf}

\begin{lemma}\label{lem2}
 For $s\in\R$ and $x_1,\dots,x_n\in X$, we have
 \begin{align*}
  \|(x_k\cos ks)|\eusm C_n\|\le\|(x_k)|\eusm C_n\|,&\quad 
  \|(x_k\sin ks)|\eusm S_n\|\le\|(x_k)|\eusm C_n\|,\\
  \|(x_k\cos ks)|\eusm S_n\|\le\|(x_k)|\eusm S_n\|,&\quad
  \|(x_k\sin ks)|\eusm C_n\|\le\|(x_k)|\eusm S_n\|.
 \end{align*}
\end{lemma}

\begin{pf}
 It follows from
 \[ 2\cos ks\cos kt=\cos k(s+t)+\cos k(s-t)
 \] and the translation invariance of the Lebesgue measure that
 \begin{align*}
  2&\left\|(x_k\cos ks)|\eusm C_n\right\|=
  \left(\frac1\pi\int_{-\pi}^{+\pi}\left\|\sum_{k=1}^nx_k\,2\cos ks
  \cos kt\right\|^2dt\right)^{1/2}\\
  &\le\left(\frac1\pi\int_{-\pi}^{+\pi}\left\|\sum_{k=1}^nx_k\cos
  k(t+s)\right\|^2dt\right)^{1/2}+\left(\frac1\pi\int_{-\pi}^{+\pi}
  \left\|\sum_{k=1}^nx_k\cos k(t-s)\right\|^2dt\right)^{1/2}\\
  &=2\|(x_k)|\eusm C_n\|.
 \end{align*}
 This proves the first inequality. The others can be proved in the same way.
\end{pf}

\begin{lemma}\label{lem3}
 $\r(\eusm S_n\otimes\eusm S_n,\eusm C_n)\le\sqrt2$.
\end{lemma}

\begin{pf}
 Squaring the inequality
 \[ \|(x_k\sin ks)|\eusm S_n\|\le\|(x_k)|\eusm C_n\|,
 \] from Lemma \ref{lem2} and integrating over $s\in[-\pi,+\pi]$ yields
 \[ \|(x_k)|\eusm S_n\otimes\eusm S_n\|^2\le2\|(x_k)|\eusm C_n\|^2.
 \] This proves the assertion.
\end{pf}

We are now ready to prove our first result.

\begin{proposition}\label{prop1}
 $\r(\eusm S_n,\eusm C_n)\le2\r(\eusm C_n,\eusm S_n)$.
\end{proposition}

\begin{pf}
 By Lemmas \ref{lem1} and \ref{lem3} as well as Proposition \ref{prop tensor},
 we get that
 \[ \r(\eusm S_n,\eusm C_n)\le
  \r(\eusm S_n,\eusm S_n\otimes\eusm C_n)\circ
  \r(\eusm S_n\otimes\eusm C_n,\eusm S_n\otimes\eusm S_n)\circ
  \r(\eusm S_n\otimes\eusm S_n,\eusm C_n)
  \le\sqrt2\,\r(\eusm C_n,\eusm S_n)\,\sqrt2.
 \] This proves the assertion.
\end{pf}

To prove the converse of Proposition \ref{prop1} we show the following
lemma.
\begin{lemma}\label{lem4}
 Let $x_1,\dots,x_n\in X$ and define $x_{-1}=x_0=x_{n+1}=x_{n+2}:=0$. Then for
 $t\in\R$, we have
 \begin{align*}
  2\sin t\sum_{k=1}^nx_k\sin kt&=\sum_{k=0}^{n+1}(x_{k+1}-x_{k-1})\cos kt,\\
  2\sin t\sum_{k=1}^nx_k\cos kt&=\sum_{k=1}^{n+1}(x_{k-1}-x_{k+1})\sin kt.
 \end{align*}
\end{lemma}

\begin{pf}
 The equations above follow from
 \begin{align*}
  2\sin t\sin kt&=\cos(k-1)t-\cos(k+1)t,\\
  2\sin t\cos kt&=\sin(k+1)t-\sin(k-1)t
 \end{align*}
 by rearranging the summation.
\end{pf}

\begin{lemma}\label{lem5}
 For $x_0,\dots,x_{n+1}\in X$, we have
 \begin{align*}
  \left(\frac2\pi\int_0^\pi\left\|\sum_{k=1}^nx_k\cos kt\right\|^2dt
  \right)^{1/2}
  &\le\left(\frac2\pi\int_0^\pi\left\|\sum_{k=0}^{n+1}x_k\cos kt
  \right\|^2dt\right)^{1/2}+\sqrt2\|x_0\|+\|x_{n+1}\|,\\
  \left(\frac2\pi\int_0^\pi\left\|\sum_{k=1}^{n+1}x_k\sin kt\right\|
  ^2dt\right)^{1/2}
  &\le\left(\frac2\pi\int_0^\pi\left\|\sum_{k=1}^nx_k\sin kt\right\|
  ^2dt\right)^{1/2}+\|x_{n+1}\|.
 \end{align*}
\end{lemma}

\begin{pf}
 We have
 \begin{align*}
  & \left( \frac2\pi \int_0^\pi \left\| \sum_{k=1}^n x_k \cos kt \right\| ^2
  dt \right) ^{1/2}\\
  \hspace*{-2pt} & \le \left( \frac2\pi \int_0^\pi \left\| \sum_{k=0}^{n+1}
  x_k \cos kt \right\| ^2 dt \right) ^{1/2} \hspace{-9pt} + \hspace{-1pt}
  \left( \frac2\pi \int_0^\pi \|x_0\| ^2 dt \right) ^{1/2} \hspace{-9pt} +
  \hspace{-1pt} \left( \frac2\pi \int_0^\pi \|x_{n+1} \cos(n+1)t\| ^2 dt
  \right) ^{1/2}\\
  \hspace*{-2pt} & = \left( \frac2\pi \int_0^\pi \left\| \sum_{k=0}^{n+1} x_k
  \cos kt \right\| ^2 dt \right) ^{1/2} \hspace{-9pt} + \hspace{-1pt} \sqrt2
  \|x_0\| + \|x_{n+1}\|.
 \end{align*}
 This proves the first inequality. The second one follows in the same way.
\end{pf}

\begin{lemma}\label{lem6}
 Let $\Delta_0:=\left[\frac\pi3,\frac{2\pi}3\right)$.
Then for $x_1,\dots,x_n\in X$
 and $T\in\boldsymbol{\frak L}(X,Y)$, we have
 \[ \left( \frac2\pi \int_{\Delta_0} \left\| \sum_{k=1}^n Tx_k \cos kt
  \right\| ^2 \right) ^{1/2} \le 4 \r(T|\eusm S_n, \eusm C_n)
  \|(x_k)|\eusm S_n\|.
 \]
\end{lemma}

\begin{pf}
 If $t\in\Delta_0$, then
 \[ \sin\tfrac\pi3\le\sin t.
 \] Moreover, by Proposition \ref{single elements} we have
 \[ \|x_1\|\le\|(x_k)|\eusm S_n\|\txt{and}\|x_n\|\le\|(x_k)|\eusm S_n\|.
 \] Applying Lemmas \ref{lem4} and \ref{lem5}, we obtain
 \begin{align*}
  &2\sin\frac\pi3\left(\frac2\pi\int_{\Delta_0}\left\|\sum_{k=1}^nTx_k\cos
  kt\right\|^2\right)^{1/2}
  \le\left(\frac2\pi\int_{\Delta_0}\left\|2\sin t\sum_{k=1}^nTx_k\cos kt
  \right\|^2\right)^{1/2}\\
  &\le\left(\frac2\pi\int_0^\pi\left\|\sum_{k=1}^{n+1}T(x_{k-1}-
  x_{k+1})\sin kt\right\|^2\right)^{1/2}\\
  &\le\left(\frac2\pi\int_0^\pi\left\|\sum_{k=1}^{n}T(x_{k-1}-x_{k+1})
  \sin kt\right\|^2\right)^{1/2}+\|Tx_n\|\\
  &\le\r(T|\eusm S_n,\eusm C_n)\!\left[\!\left(\frac2\pi\int_0^\pi
  \left\|\sum_{k=1}^n(x_{k-1}-x_{k+1})\cos kt\right\|^2\!dt\right)
  ^{1/2}+\|x_n\|\right]\\
  &\le\r(T|\eusm S_n,\eusm C_n)\!\left[\!\left(\frac2\pi\int_0^\pi
  \left\|\sum_{k=0}^{n+1}(x_{k-1}-x_{k+1})\cos kt\right\|^2\!dt\right)
  ^{1/2}+\sqrt2\|x_1\|+2\|x_n\|\right]\\
  &=\r(T|\eusm S_n,\eusm C_n)\!\left[2\left(\frac2\pi\int_0^\pi
  \left\|\sin t\sum_{k=1}^nx_k\sin kt\right\|^2\!dt\right)^{1/2}+
  \sqrt2\|x_1\|+2\|x_n\|\right]\\
  &\le\r(T|\eusm S_n,\eusm C_n)\left[2\|(x_k)|\eusm S_n\|+\sqrt2\|x_1\|
  +2\|x_n\|\right]\le(4+\sqrt2)\r(T|\eusm S_n,\eusm C_n)
  \|(x_k)|\eusm S_n\|.
 \end{align*}
 The preceding inequality yields the assertion, since
 \[ \frac{4+\sqrt2}{2\sin\frac\pi3}=3.1258\ldots<4.\qed
 \]
 \renewcommand{\qed}{}
\end{pf}

\begin{proposition}\label{prop2}
 $\r(\eusm C_n,\eusm S_n)\le9\r(\eusm S_n,\eusm C_n)$.
\end{proposition}

\begin{pf}
 Obviously,
 \begin{equation}\label{eq1}
  \|(Tx_k)|\eusm C_n\|=\left(\frac2\pi\int_0^\pi\left\|\sum_{k=1}^nTx_k
  \cos kt\right\|^2dt\right)^{1/2}=\left(I_{-1}^2+I_0^2+I_{+1}^2\right)
  ^{1/2},
 \end{equation}
 where
 \[ I_\alpha:=\left(\frac2\pi\int_{\Delta_\alpha}\left\|\sum_{k=1}^nTx_k\cos
  kt\right\|^2dt\right)^{1/2}\txt{and}\Delta_\alpha
  :=\left[\frac13\pi,\frac23\pi\right)+\frac\alpha3\pi.
 \] We know from Lemma \ref{lem6} that
 \begin{equation}\label{eq2}
  I_0\le4\r(T|\eusm S_n,\eusm C_n)\|(x_k)|\eusm S_n\|.
 \end{equation}
 In order to estimate $I_\alpha$ with $\alpha=\pm1$, we
substitute $s:=t\mp\frac\pi3$. Then
 \begin{align*}
  I_\alpha&=\left(\frac2\pi\int_{\Delta_0}\left\|\sum_{k=1}^nTx_k(\cos
  k\tfrac\pi3\cos ks\mp\sin k\tfrac\pi3\sin ks)\right\|^2ds
  \right)^{1/2}\\
  &\le\left(\frac2\pi\int_{\Delta_0}\left\|\sum_{k=1}^nT(x_k\cos
  k\tfrac\pi3)\cos ks\right\|^2ds\right)^{1/2}\hspace{-5pt}+
  \left(\frac2\pi\int_{\Delta_0}\left\|\sum_{k=1}^nT(x_k\sin
  k\tfrac\pi3)\sin ks\right\|^2ds\right)^{1/2}\hspace{-3pt}.
 \end{align*}
 We now estimate the first summand by Lemma \ref{lem6} and the second summand 
 by applying the defining inequality \eqref{def r} of $\r(T|\eusm S_n,\eusm 
 C_n)$. This yields
 \[ I_\alpha\le\r(T|\eusm S_n,\eusm C_n)\left(4\left\|\left.(x_k\cos 
  k\tfrac\pi3)\right|\eusm S_n\right\|+\left\|\left.(x_k\sin 
  k\tfrac\pi3)\right|\eusm C_n\right\|\right).
 \] Hence, in view of Lemma \ref{lem2}, we arrive at
 \begin{equation}\label{eq3}
  I_\alpha\le5\,\r(T|\eusm S_n,\eusm C_n)\|(x_k)|\eusm S_n\|.
 \end{equation}
 Combining \eqref{eq1}, \eqref{eq2} and \eqref{eq3} yields
 \[ \|(Tx_k)|\eusm C_n\|\le(25+16+25)^{1/2}\r(T|\eusm S_n,\eusm 
  C_n)\|(x_k)|\eusm S_n\|.
 \] In view of
 \[ \sqrt{66}=8.1240\ldots<9,
 \] this completes the proof.
\end{pf}

\begin{remark}
 By using the exact value of $3.1258\ldots$ for the constant appearing in 
 Lemma \ref{lem6}, we can even obtain a value of $6.6194\ldots$ for the 
 constant in the previous proposition.
\end{remark}

We now deal with the ideal norms $\d(\eusm E_n,\eusm E_n)$.

\begin{lemma}\label{lem mono}
 For $m,n\in\N$ with $n<m$, we have
 \[ \d(\eusm E_{m\pm n},\eusm E_{m\pm n})\le\d(\eusm E_m,\eusm E_m)+
  \d(\eusm E_n,\eusm E_n).
 \]
\end{lemma}

\begin{pf}
 We have
 \begin{align*}
  \|(T\sprod f{\overline e_k})|\eusm E_{m+n}\|
  &\le\|(T\sprod f{\overline e_k})|\eusm E_m\|+
  \|(T\sprod f{\overline e_{m+k}})|\eusm E_n\|\\
  &=\|(T\sprod f{\overline e_k})|\eusm E_m\|+
 \|(T\sprod{fe_{-m}}{\overline e_k})|\eusm E_n\|\\
  &\le\d(T|\eusm E_m,\eusm E_m)\,\|f|L_2\|+
  \d(T|\eusm E_n,\eusm E_n)\,\|fe_{-m}|L_2\|\\
  &\le\big(\d(T|\eusm E_m,\eusm E_m)+\d(T|\eusm E_n,\eusm E_n)\big)
  \|f|L_2\|.\\
 \intertext{Similarly,}
  \|(T\sprod f{\overline e_k})|\eusm E_{m-n}\|
  &\le\|(T\sprod f{\overline e_k})|\eusm E_m\|+
  \|(T\sprod f{\overline e_{m-n+k}})|\eusm E_n\|\\
  &=\|(T\sprod f{\overline e_k})|\eusm E_m\|+
  \|(T\sprod{fe_{-m+n}}{\overline e_k})|\eusm E_n\|\\
  &\le\d(T|\eusm E_m,\eusm E_m)\,\|f|L_2\|+
  \d(T|\eusm E_n,\eusm E_n)\,\|fe_{-m+n}|L_2\|\\
  &\le\big(\d(T|\eusm E_m,\eusm E_m)+\d(T|\eusm E_n,\eusm E_n)\big)
  \|f|L_2\|.
 \end{align*}
\end{pf}

In the following, we write
\[ \m_n(T) := \max \big\{ \r(T| \eusm C_n, \eusm S_n), \r(T| \eusm S_n,
 \eusm C_n) \big\}.
\]

\begin{lemma}\label{lem7}
 For $x_{-n},\dots,x_0,\dots,x_{+n}\in X$, we have
 \[ \left(\frac1{2\pi}\int_{-\pi}^{+\pi}\left\|\sum_{k=1}^nTx_k\exp(ikt)
  \right\|^2dt\right)^{1/2}\le 4\,\m_n(T)\left(\frac1{2\pi}\int_{-\pi}
  ^{+\pi}\left\|\sum_{|k|\le n}x_k\exp(ikt)\right\|^2dt
  \right)^{1/2}.
 \]
\end{lemma}

\begin{pf}
 For $k=1,\dots,n$, we let
 \[ u_k:=x_k+x_{-k}\txt{and}v_k:=x_k-x_{-k}.
 \] It follows from $u_k+v_k=2x_k$ and Euler's formula that
 \[ x_k\exp(ikt)=\tfrac{1}{2}\big(u_k\cos kt+v_k\cos kt+iu_k\sin kt+iv_k
  \sin kt\big).
 \] Hence
 \begin{align*}
  \begin{split}
   &\|(Tx_k)|\eusm E_n\|=\left(\frac1{2\pi}\int_{-\pi}^{+\pi}
   \left\|\sum_{k=1}^nTx_k\exp(ikt)\right\|^2dt\right)^{1/2}\\[5pt]
   &\le\frac12\left\{
   \begin{array}{clcl}
    &\left(\frac1{2\pi}\int\limits_{-\pi}^{+\pi}\|\sum\limits_{k=1}^nTu_k
     \cos kt\|^2dt \right)^{1/2}&
    +&\left(\frac1{2\pi}\int\limits_{-\pi}^{+\pi}\|\sum\limits_{k=1}^nTv_k
     \cos kt\|^2dt \right)^{1/2}\\
    +&\left(\frac1{2\pi}\int\limits_{-\pi}^{+\pi}\|\sum\limits_{k=1}^nTu_k
     \sin kt\|^2dt \right)^{1/2}&
    +&\left(\frac1{2\pi}\int\limits_{-\pi}^{+\pi}\|\sum\limits_{k=1}^nTv_k
     \sin kt\|^2dt \right)^{1/2}
   \end{array}
   \right\}\displaybreak[0]\\[5pt]
   &\le\frac12\left\{
   \begin{array}{ccl}
    &\|T\|&\left(\frac1{2\pi}\int\limits_{-\pi}^{+\pi}\|\sum\limits_{k=1}^nu_k
     \cos kt\| ^2dt\right)^{1/2}\\
    +&\r(T|\eusm C_n,\eusm S_n)&\left(\frac1{2\pi}\int\limits_{-\pi}^{+\pi}\|
     \sum\limits_{k=1}^nv_k\sin kt\|^2dt\right)^{1/2}\\
    +&\r(T|\eusm S_n,\eusm C_n)&\left(\frac1{2\pi}\int\limits_{-\pi}^{+\pi}\|
     \sum\limits_{k=1}^nu_k\cos kt\|^2dt\right)^{1/2}\\
    +&\|T\|&\left(\frac1{2\pi}\int\limits_{-\pi}^{+\pi}\|\sum\limits_{k=1}^nv_k
     \sin kt\| ^2dt\right)^{1/2}
   \end{array}
   \right\}\\[5pt]
   &\le\m_n(T)\left[\left(\tfrac1{2\pi}\int_{-\pi}^{+\pi}\left\|
    \sum_{k=1}^nu_k\cos kt\right\|^2dt\right)^{1/2}
   +\left(\tfrac1{2\pi}\int_{-\pi}^{+\pi}\left\|\sum_{k=1}^nv_k\sin kt
    \right\|^2dt\right)^{1/2}\right].
  \end{split}
 \end{align*}
 By the obvious fact that
 \[ \|u\|,\|v\| \le \frac12 \left( \|u+iv\|+\|u-iv\| \right) 
  \txt{for $u,v\in X$,}
 \] we obtain
 \[ \|(Tx_k)|\eusm E_n\|\le \m_n(T)\left\{
  \begin{array}{c}
   \left(\frac1{2\pi}\int\limits_{-\pi}^{+\pi}\|\sum\limits_{k=1}^n(u_k\cos 
    kt +iv_k\sin kt)\|^2dt\right)^{1/2}\\+\\
   \left(\frac1{2\pi}\int\limits_{-\pi}^{+\pi}\|\sum\limits_{k=1}^n(u_k\cos 
    kt -iv_k\sin kt)\|^2dt\right)^{1/2}
  \end{array}
  \right\}.
 \] Substituting $-t$ for $t$ in the lower term on the right--hand side
 yields
 \[ \|(Tx_k)|\eusm E_n\|\le2\,\m_n(T)\left(\frac1{2\pi}\int_{-\pi}^{+\pi}
  \left\|\sum_{k=1}^n(u_k\cos kt+iv_k\sin kt)\right\|^2dt\right)^{1/2}.
 \] Finally, we conclude from
 \[ \sum_{|k|\le n}x_k\exp(ikt)-x_0=\sum_{k=1}^n(u_k\cos kt+iv_k\sin kt)
 \] and Proposition \ref{single elements} that
 \begin{align*}
  \left(\frac1{2\pi}\int_{-\pi}^{+\pi}\left\|
  \sum_{k=1}^n(u_k\cos kt+iv_k\sin kt)
  \right\|^2dt\right)^{1/2}\!\!&\le
  \left(\frac1{2\pi}\int_{-\pi}^{+\pi}\left\|\sum_{|k|\le 
  n}x_k\exp(ikt)\right\|^2
  dt\right)^{1/2}\!\!+\|x_0\|\\
  &\le2\left(\frac1{2\pi}\int_{-\pi}^{+\pi}\left\|\sum_{|k|\le 
  n}x_k\exp(ikt)\right\|^2
  dt\right)^{1/2}.
 \end{align*}
 This proves the desired result.
\end{pf}

The basic trick in the following proof goes back to M.~Junge.

\begin{proposition}\label{prop3}
 $\d(\eusm E_n,\eusm E_n)\le96\,\m_n$.
\end{proposition}

\begin{pf}
 The $m$--th de la Vall\'ee Poussin kernel $V_m$ is defined as
 \[ V_m(t):=\frac1m\sum_{k=m}^{2m-1}D_k(t),
 \] where
 \[ D_k(t):=\sum_{|l|\le k}\exp(ilt)
 \] is the $k$--th Dirichlet kernel. It is known that
 \begin{equation}\label{vallee koeff}
  \sprod{V_m}{\overline e_k}=\frac1{2\pi}\int_{-\pi}^{+\pi}V_m(t)
  \overline{e_k(t)}\,dt
  =\left\{
   \begin{array}{cl}
    1&\txt{for $|k|\le m$,}\\[2pt]
    \frac{2m-|k|}m&\txt{for $m<|k|<2m$,}\\[2pt]
    0&\txt{for $|k|\ge2m$,}
   \end{array}
   \right.
 \end{equation}
 and
 \begin{equation}\label{vallee norm}
  \frac1{2\pi}\int_{-\pi}^{+\pi}|V_m(t)|\,dt\le3;
 \end{equation}
 see e.g.~Zygmund \cite{ZYG}.

 On $L_2^X(-\pi,+\pi)$ we consider the $m$--th de la Vall\'ee Poussin operator
 \[ V_m:f(t)\longrightarrow\frac1{2\pi}\int_{-\pi}^{+\pi}V_m(s-t)f(t)\,dt.
 \] It follows from \eqref{vallee norm} that
 \[ \|V_mf|L_2\|\le3\,\|f|L_2\|.
 \] For $f\in L_2^X(-\pi,+\pi)$, we let
 \[ x_k^{(m)}:=\sprod{V_mf}{\overline e_k}=\sprod{V_m}{\overline e_k}\sprod f
  {\overline e_k}.
 \] Hence, by \eqref{vallee koeff} 
 \[ V_m f = \sum_{|k|\le2m-1} e_k \otimes x_k^{(m)} .
 \] The triangle inequality implies that
 \[ \| (T \sprod f{\overline e_k}) |\eusm E_m\| = \left( \frac1{2\pi}
  \int_{-\pi}^{+\pi} \left\| \sum_{k=1}^m Tx_k^{(m)} e_k(t) \right\| ^2 dt
  \right) ^{1/2} \le I_1 + I_2,
 \] where
 \[ I_1:=\left(\frac1{2\pi}\int_{-\pi}^{+\pi}
  \left\|\sum_{k=1}^{2m-1}Tx_k^{(m)}e_k(t)\right\|^2dt\right)^{1/2}
  \text{and}\,\,I_2:=\left(\frac1{2\pi}\int_{-\pi}^{+\pi}
  \left\|\sum_{k=m+1}^{2m-1}Tx_k^{(m)}e_k(t)\right\|^2dt\right)^{1/2}
  \hspace{-2pt}.
 \] Lemma \ref{lem7} implies that
 \begin{align*}
  I_1&\le4\,\m_{2m-1}(T)\left(\frac1{2\pi}\int_{-\pi}^{+\pi}\left\|\sum
  _{|k|\le2m-1}x_k^{(m)}\exp(ikt)\right\|^2dt\right)^{1/2}\\
  &=4\,\m_{2m-1}(T)\|V_mf|L_2\|\le12\,\m_{2m-1}(T)\|f|L_2\|.
 \end{align*}
 To estimate the second term, we recall that $x_k^{(m)}=0$ if $|k|\ge2m$.
 Therefore
 \begin{align*}
  I_2&=\left(\frac1{2\pi}\int_{-\pi}^{+\pi}\left\|\sum_{k=m+1}^{2m-1}
  Tx_k^{(m)}\exp(ikt)\right\|^2dt\right)^{1/2}\\
  &=\left(\frac1{2\pi}\int_{-\pi}^{+\pi}\left\|\sum_{k=1}^{3m-1}
  Tx_{k+m}^{(m)}\exp(ikt)\right\|^2dt\right)^{1/2}\\
  &\le4\,\m_{3m-1}(T)\left(\frac1{2\pi}\int_{-\pi}^{+\pi}\left\|\sum
  _{|k|\le3m-1}x_{k+m}^{(m)}\exp(ikt)\right\|^2dt\right)^{1/2}\\
  &=4\,\m_{3m-1}(T)\left(\frac1{2\pi}\int_{-\pi}^{+\pi}\left\|\sum
  _{|k|\le2m-1}x_k^{(m)}\exp(ikt)\right\|^2dt\right)^{1/2}\\
  &=4\,\m_{3m-1}(T)\|V_mf|L_2\|\le12\,\m_{3m-1}(T)\|f|L_2\|.
 \end{align*}
 Combining the preceding estimates and taking into account the monotonicity of
 $\m_n(T)$, we arrive at
 \[ \d(T|\eusm E_m,\eusm E_m)\le24\,\m_{3m-1}(T).
 \] To complete the proof, for given $n\in\N$, we choose $m$ such that
 $3m-1\le n\le 3m+1$. Then it follows from Lemma \ref{lem mono} that
 \[ \d(\eusm E_n,\eusm E_n)=\left\{
  \begin{array}{lcl}
   \d(\eusm E_{3m-1},\eusm E_{3m-1})&\le&2\,\d(\eusm E_m,\eusm E_m)+
\d(\eusm E_{m-1},\eusm E_{m-1})\\
   \d(\eusm E_{3m},\eusm E_{3m})&\le&3\,\d(\eusm E_m,\eusm E_m)\\
   \d(\eusm E_{3m+1},\eusm E_{3m+1})&\le&2\,
\d(\eusm E_m,\eusm E_m)+\d(\eusm E_{m+1},\eusm E_{m+1})
  \end{array}
  \right\}\le4\,\d(\eusm E_m,\eusm E_m).
 \] Hence
 \[ \d(\eusm E_n,\eusm E_n)\le4\,\d(\eusm E_m,\eusm E_m)\le96\,\m_{3m-1}\le
  96\,\m_n.
 \]
\end{pf}

The next proposition is a special case of \eqref{r le d}.
\begin{proposition}\label{prop4}
 $\r(\eusm S_n,\eusm C_n)\le\d(\eusm S_n,\eusm C_n)$ and
 $\r(\eusm C_n,\eusm S_n)\le\d(\eusm C_n,\eusm S_n)$.
\end{proposition}

To estimate the ideal norms $\d(\eusm S_n,\eusm C_n)$ and $\d(\eusm 
C_n,\eusm S_n)$ by $\d(\eusm E_n,\eusm E_n)$, we need one more lemma.

\begin{lemma}\label{lem8}
 $\r(\eusm C_n,\eusm E_n)\le\sqrt2$ and $\r(\eusm S_n,\eusm E_n)\le\sqrt2$.
\end{lemma}

\begin{pf}
 By Euler's formula, we have
 \[ c_k=\tfrac1{\sqrt2}(e_k+\overline e_k).
 \] Hence by \eqref{conju}
 \[ \|(x_k)|\eusm C_n\|\le\tfrac1{\sqrt2}\left(\|(x_k)|\eusm E_n\|+
  \|(x_k)|\overline{\eusm E}_n\|\right)=\sqrt2\;\|(x_k)|\eusm E_n\|.
 \] This proves the left--hand inequality. The right--hand inequality can be 
 obtained in the same way.
\end{pf}

The next proposition follows immediately from Proposition \ref{cor} and 
Lemma \ref{lem8}.
\begin{proposition}\label{prop5}
 $\d(\eusm S_n,\eusm C_n)\le2\d(\eusm E_n,\eusm E_n)$ and
 $\d(\eusm C_n,\eusm S_n)\le2\d(\eusm E_n,\eusm E_n)$.
\end{proposition}

We now combine Propositions \ref{prop1} through \ref{prop5} to complete the 
proof of the theorem.
\begin{pf*}{Proof of the theorem}
 Proposition \ref{prop3} states that $\d(\eusm E_n,\eusm E_n)$ lies below
 $\r(\eusm S_n,\eusm C_n)$ and $\r(\eusm C_n,\eusm S_n)$.

 Proposition \ref{prop4} implies that the sequences $\r(\eusm S_n,\eusm C_n)$ 
 and $\r(\eusm C_n,\eusm S_n)$, respectively, lie below $\d(\eusm S_n,\eusm 
 C_n)$ and $\d(\eusm C_n,\eusm S_n)$.

 Finally, it follows from Proposition \ref{prop5} that $\d(\eusm S_n,\eusm 
 C_n)$ and $\d(\eusm C_n,\eusm S_n)$ lie below $\d(\eusm E_n,\eusm E_n)$.

 This proves the uniform equivalence of all five sequences of ideal norms and 
 thus completes the proof of the theorem.
\end{pf*}


\ifx\undefined\bysame
\newcommand{\bysame}{\leavevmode\hbox to3em{\hrulefill}\,}
\fi



\begin{thebibliography}{1}

\bibitem{Bou83}
Bourgain, J.: {\em Some remarks on {B}anach spaces in which martingale
  difference sequences are unconditional}, Ark. Mat. {\bf 22} (1983), 163--168.

\bibitem{BUR:1986}
Burkholder, D.~L.: {\em Martingales and {F}ourier analysis in {B}anach
  spaces}, Probability and analysis (Varenna, Italy, 1985), Lecture Notes in
  Mathematics, no. 1206, Springer--Verlag, 1986, pp.~61--108.

\bibitem{DEF:1986}
Defant, M.: {\em Zur vektorwertigen {H}ilberttransformation}, Ph.D. thesis,
  Christian--Albrechts--Universit\"at Kiel, 1986.

\bibitem{DIN}
Dinculeanu, N.: {\em Vector measures}, Deutscher Verlag der Wissenschaften,
  Berlin, 1966.

\bibitem{PIE:1}
Pietsch, A.: {\em Operator ideals}, Deutscher Verlag der Wissenschaften,
  Berlin, 1978.

\bibitem{PIE:2}
\bysame, {\em Eigenvalues and $s$--numbers}, Akademische Verlagsgesellschaft
  Geest \& Portig K.--G., Leipzig, 1987.

\bibitem{PIW}
Pietsch, A. and Wenzel, J.: {\em Orthogonal sytems and geometry of
  {B}anach spaces}, In preparation.

\bibitem{ZYG}
Zygmund, A.: {\em Trigonometric series}, second ed., Cambridge University 
Press, Cambridge, 1959.

\end{thebibliography}
\end{document}